\documentclass[10pt,a4paper]{article}
\usepackage{amssymb}
\usepackage{amscd}
\usepackage{latexsym}
\usepackage{amsmath}
\usepackage{amsfonts}
\usepackage{amscd}

\begin{document}

\newtheorem{theorem}{Theorem}[section]
\newtheorem{remark}[theorem]{Remark}
\newtheorem{mtheorem}[theorem]{Main Theorem}
\newtheorem{bbtheo}[theorem]{The Strong Black Box}
\newtheorem{observation}[theorem]{Observation}
\newtheorem{proposition}[theorem]{Proposition}
\newtheorem{lemma}[theorem]{Lemma}
\newtheorem{testlemma}[theorem]{Test Lemma}
\newtheorem{mlemma}[theorem]{Main Lemma}
\newtheorem{note}[theorem]{{\bf Note}}
\newtheorem{steplemma}[theorem]{Step Lemma}
\newtheorem{corollary}[theorem]{Corollary}
\newtheorem{notation}[theorem]{Notation}
\newtheorem{example}[theorem]{Example}
\newtheorem{definition}[theorem]{Definition}

\renewcommand{\labelenumi}{(\roman{enumi})}
\def\Pf{\smallskip\goodbreak{\sl Proof. }}

\def\Fin{\mathop{\rm Fin}\nolimits}
\def\br{\mathop{\rm br}\nolimits}
\def\fin{\mathop{\rm fin}\nolimits}
\def\Ann{\mathop{\rm Ann}\nolimits}
\def\Aut{\mathop{\rm Aut}\nolimits}
\def\End{\mathop{\rm End}\nolimits}
\def\bfb{\mathop{\rm\bf b}\nolimits}
\def\bfi{\mathop{\rm\bf i}\nolimits}
\def\bfj{\mathop{\rm\bf j}\nolimits}
\def\df{{\rm df}}
\def\bfk{\mathop{\rm\bf k}\nolimits}
\def\bEnd{\mathop{\rm\bf End}\nolimits}
\def\iso{\mathop{\rm Iso}\nolimits}
\def\id{\mathop{\rm id}\nolimits}
\def\Ext{\mathop{\rm Ext}\nolimits}
\def\Ines{\mathop{\rm Ines}\nolimits}
\def\Hom{\mathop{\rm Hom}\nolimits}
\def\bHom{\mathop{\rm\bf Hom}\nolimits}
\def\Rk{ R_\k-\mathop{\bf Mod}}
\def\Rn{ R_n-\mathop{\bf Mod}}
\def\map{\mathop{\rm map}\nolimits}
\def\cf{\mathop{\rm cf}\nolimits}
\def\top{\mathop{\rm top}\nolimits}
\def\Ker{\mathop{\rm Ker}\nolimits}
\def\Bext{\mathop{\rm Bext}\nolimits}
\def\Br{\mathop{\rm Br}\nolimits}
\def\dom{\mathop{\rm Dom}\nolimits}
\def\min{\mathop{\rm min}\nolimits}
\def\im{\mathop{\rm Im}\nolimits}
\def\max{\mathop{\rm max}\nolimits}
\def\rk{\mathop{\rm rk}}
\def\Diam{\diamondsuit}
\def\Z{{\mathbb Z}}
\def\Q{{\mathbb Q}}
\def\N{{\mathbb N}}
\def\bQ{{\bf Q}}
\def\bF{{\bf F}}
\def\bX{{\bf X}}
\def\bY{{\bf Y}}
\def\bHom{{\bf Hom}}
\def\bEnd{{\bf End}}
\def\bS{{\mathbb S}}
\def\AA{{\cal A}}
\def\BB{{\cal B}}
\def\CC{{\cal C}}
\def\DD{{\cal D}}
\def\TT{{\cal T}}
\def\FF{{\cal F}}
\def\GG{{\cal G}}
\def\PP{{\cal P}}
\def\SS{{\cal S}}
\def\XX{{\cal X}}
\def\YY{{\cal Y}}
\def\fS{{\mathfrak S}}
\def\fH{{\mathfrak H}}
\def\fU{{\mathfrak U}}
\def\fW{{\mathfrak W}}
\def\fK{{\mathfrak K}}
\def\PT{{\mathfrak{PT}}}
\def\T{{\mathfrak{T}}}
\def\fX{{\mathfrak X}}
\def\fP{{\mathfrak P}}
\def\X{{\mathfrak X}}
\def\Y{{\mathfrak Y}}
\def\F{{\mathfrak F}}
\def\C{{\mathfrak C}}
\def\B{{\mathfrak B}}
\def\J{{\mathfrak J}}
\def\fN{{\mathfrak N}}
\def\fM{{\mathfrak M}}
\def\Fk{{\F_\k}}
\def\bar{\overline }
\def\Bbar{\bar B}
\def\Cbar{\bar C}
\def\Pbar{\bar P}
\def\etabar{\bar \eta}
\def\Tbar{\bar T}
\def\fbar{\bar f}
\def\nubar{\bar \nu}
\def\rhobar{\bar \rho}
\def\Abar{\bar A}
\def\a{\alpha}
\def\b{\beta}
\def\g{\gamma}
\def\w{\omega}
\def\e{\varepsilon}
\def\o{\omega}
\def\va{\varphi}
\def\k{\kappa}
\def\m{\mu}
\def\n{\nu}
\def\r{\rho}
\def\f{\phi}
\def\hv{\widehat\v}
\def\hF{\widehat F}
\def\v{\varphi}
\def\s{\sigma}
\def\l{\lambda}
\def\lo{\lambda^{\aln}}
\def\d{\delta}
\def\z{\zeta}
\def\th{\theta}
\def\t{\tau}
\def\ale{\aleph_1}
\def\aln{\aleph_0}
\def\Cont{2^{\aln}}
\def\nld{{}^{ n \downarrow }\l}
\def\n+1d{{}^{ n+1 \downarrow }\l}
\def\hsupp#1{[[\,#1\,]]}
\def\size#1{\left|\,#1\,\right|}
\def\Binfhat{\widehat {B_{\infty}}}
\def\Zhat{\widehat \Z}
\def\Mhat{\widehat M}
\def\Rhat{\widehat R}
\def\Phat{\widehat P}
\def\Fhat{\widehat F}
\def\fhat{\widehat f}
\def\Ahat{\widehat A}
\def\Chat{\widehat C}
\def\Ghat{\widehat G}
\def\Bhat{\widehat B}
\def\Btilde{\widetilde B}
\def\Ftilde{\widetilde F}
\def\restr{\mathop{\upharpoonright}}
\def\to{\rightarrow}
\def\arr{\longrightarrow}
\def\LA{\langle}
\def\RA{\rangle}
\newcommand{\norm}[1]{\text{$\parallel\! #1 \!\parallel$}}
\newcommand{\supp}[1]{\text{$\left[ \, #1\, \right]$}}
\def\set#1{\left\{\,#1\,\right\}}
\newcommand{\mb}{\mathbf}
\newcommand{\wt}{\widetilde}
\newcommand{\card}[1]{\mbox{$\left| #1 \right|$}}
\newcommand{\union}{\bigcup}
\newcommand{\inters}{\bigcap}
\newcommand{\ER}{{\rm E}}
\def\Proof{{\sl Proof.}\quad}
\def\fine{\ \black\vskip.4truecm}
\def\black{\ {\hbox{\vrule width 4pt height 4pt depth
0pt}}}
\def\fine{\ \black\vskip.4truecm}
\long\def\alert#1{\smallskip\line{\hskip\parindent\vrule%
\vbox{\advance\hsize-2\parindent\hrule\smallskip\parindent.4\parindent%
\narrower\noindent#1\smallskip\hrule}\vrule\hfill}\smallskip}

\title{Solution to the Uniformly Fully Inert \\ Subgroups Problem for Abelian Groups}
\author{Andrey R. Chekhlov \\Faculty of Mathematics and Mechanics, Section of Algebra, \\Tomsk State University, Tomsk 634050, Russia\\{\small e-mails: cheklov@math.tsu.ru, a.r.che@yandex.ru}
\\and\\ Peter V. Danchev \\Institute of Mathematics and Informatics, Section of Algebra, \\Bulgarian Academy of Sciences, Sofia 1113, Bulgaria\\{\small e-mails: danchev@math.bas.bg, pvdanchev@yahoo.com}}
\maketitle

\centerline{(To Brendan Goldsmith on his {\bf 75}th birthday)}

\medskip
\medskip

\begin{abstract}{A famous conjecture attributed to Dardano-Dikranjan-Rinauro-Salce states that {\it any uniformly fully inert subgroup of a given group is commensurable with a fully invariant subgroup} (see, respectively, \cite{DDR} and \cite{DDS}). In this short note, we completely settle this problem in the affirmative for an arbitrary Abelian group.}
\end{abstract}

\section{Preliminaries}

Throughout the present brief paper, all our groups are {\it additively} written and {\it Abelian}. Our notation and terminology are mainly standard and follow those from \cite{Kap}.

\medskip

Recall the standard notion that a subgroup $F$ of an arbitrary group $G$ is said to be {\it fully invariant} if $\phi(F)\subseteq F$ for any endomorphism $\phi$ of $G$, while if $\phi$ is an invertible endomorphism (= an automorphism), then $F$ is said to be a {\it characteristic} subgroup. It is obvious that fully invariant subgroups are always characteristic, whereas the converse implication fails in general.

\medskip

On the same vein, imitating \cite{DDS} and \cite{CDG}, respectively, a subgroup $S$ of a group $G$ is called {\it fully (resp., characteristically) inert}, provided $(\phi(S)+S)/S$ is finite for all endomorphisms, respectively automorphisms, $\phi$ of $G$. 

The first of these two concepts is known to be refined in \cite{DDR} and \cite{DDS} to the so-called {\it uniformly fully inert} subgroups by requiring the existence of a fixed positive integer $m$ such that the cardinality of the quotient-group $(\phi(S)+S)/S$ is bounded by $m$ (i.e., it has at most $m$ elements, for each endomorphism $\phi$ of $G$), while the second concept mentioned above was refined in \cite{DDR} by defining a subgroup $H$ of a group $G$ to be {\it uniformly characteristically inert}, provided there is a positive integer $k$ such that, for every automorphism $\varphi$ of $G$, the cardinality of the factor-group $(\varphi(H)+H)/H$ is bounded by $k$ (that is, it has no more than $k$ elements for each automorphism $\varphi$ of $G$). It is evident that any (uniformly) fully inert subgroup is always (uniformly) characteristically inert, whereas the reverse implication is wrong in general.

\medskip

Remember also that two subgroups $B$ and $C$ of a group $G$ are called {\it commensurable} and write for short that $B\sim C$ or, equivalently, that $C\sim B$ since this is a symmetric relation, provided that both quotients $(B+C)/B$ and $(B+C)/C$ are finite.

\medskip

The class ${\mathcal I}$ of groups in which fully inert subgroups are commensurable with fully invariant subgroups was studied by many authors: in fact, it is worthwhile noticing that, in many special cases, a fully inert subgroup is commensurable with a fully invariant subgroups as well as a characteristically inert subgroup is commensurable with a characteristic subgroup (see, e.g., \cite{GS1}, \cite{GS2}, \cite{GSZ} and, respectively, \cite{CDGM}, \cite{CDG}, \cite{DK}). Note that some of the results from \cite{CDGM} and \cite{CDG} that are related to the current subject, namely \cite[Proposition 4.2]{CDG} and \cite[Lemma 2.1]{CDGM}, are formulated only for $p$-groups, where $p$ is a fixed prime, but actually they remain true for arbitrary groups.

\medskip

However, it is not so hard to construct examples such that a fully inert subgroup is {\it not} commensurable with a fully invariant subgroup as well as a characteristically inert subgroup that is {\it not} commensurable with a characteristic subgroup (see, for instance, the previously cited four articles).

\medskip

On the other side, it was shown in \cite[Corollary 1.9]{DDR} the important fact that {\it any uniformly characteristically inert subgroup is commensurable a characteristic subgroup} and, moreover, in a way of similarity it was conjectured in \cite[Conjecture 5.2]{DDR} that {\it every uniformly fully inert subgroup of a given (not necessarily Abelian) group is commensurable with a fully invariant subgroup} (see \cite[Conjecture 1.6]{DDS} too).

\medskip

Our objective in the present short article is to give a complete positive solution of this difficult question in the commutative case that we illustrate in the next section.

\section{The Solution}

To simplify the exposition, let us define $\mathcal{U}$ to be the class of groups in which {\it uniformly fully inert subgroups are commensurable with fully invariant subgroups}. This class is, obviously, a counterpart of the aforementioned class $\mathcal{I}$.

\medskip

Our key tools, necessary to resolve the listed above conjecture, are the following three useful statements. The first ingredient of the proof of the main theorem listed below is the following one.

\begin{lemma}\label{1} Let $G$ be a group. Then every characteristic subgroup of $G\oplus G$ is fully invariant.
\end{lemma}

\Pf Assume $G$ is non-zero. An application of \cite[Theorem 3]{H} to the endomorphism ring $\mathrm{End}(G)$ of the group $G$ for $n=2$ guarantees that each element of $\mathrm{End}(G\oplus G)$ is the sum of three units, which means that every endomorphism of $G\oplus G$ is the sum of three automorphisms, thus immediately giving the claim.
\fine

Notice that this fact is a simple consequence of a result due to Kaplansky (see \cite[Lemma 2]{H}) pertaining to some special matrix presentations (we emphasize that this trick was also used in \cite{CDG}, \cite{DK}).

\medskip

We now state and prove the following quite curious statement.

\begin{proposition}\label{square} Every group of the form $G\oplus G$ belongs to $\mathcal{U}$.
\end{proposition}

\Pf Suppose $K := G\oplus G$ and assume $C$ is a uniformly fully inert subgroup of $K$. It is now routinely checked that $C$ is uniformly characteristically inert in $K$ (as each automorphism is necessarily an endomorphism). Therefore, with the aid of \cite[Corollary 1.9]{DDR}, we deduce that $C$ is commensurable with a characteristic subgroup $H$ of $K$. However, with Lemma~\ref{1} in hand, it then follows that $H$ is fully invariant in $K$, as needed.
\fine

Our next pivotal instrument is the following one.

\begin{proposition}\label{2} If $A$ is a group such that $A\oplus A\in {\mathcal{U}}$, then $A\in \mathcal{U}$.
\end{proposition}

\Pf The idea is based on these two well-known facts established in \cite{DDR} and \cite{DDS}, respectively (see too the literature given there):

\medskip

\noindent {\bf Fact 1:} If $H$ is a uniformly fully inert subgroup of $G$, then $H\oplus  H$ is uniformly fully inert in $G\oplus  G$.

\medskip

\noindent{\bf Fact 2:} If $A$ and $B$ are subgroups of $G$ such that $A\oplus A$ is commensurable with $B\oplus B$ in $G\oplus G$, then $A$ is commensurable with $B$ in $G$.

\medskip

With these two claims at hand, it is now straightforwardly verified that the statement is true, as promised. In fact, suppose $U$ is a uniformly fully inert subgroup of $A$. By Fact 1, the square $U\oplus U$ is uniformly fully inert subgroup of $A\oplus A$. So, by assumption, $U\oplus U$ is commensurable with a fully invariant subgroup of $A\oplus A$, say $V\oplus V$ for some $V\leq A$. Thus, Fact 2 is now applicable to get that $U$ is commensurable with $V$ in $A$, where it is readily checked that $V$ is fully invariant in $A$, as required.
\fine

We now come to our basic assertion that answers in the affirmative the important conjecture from respectively \cite[Conjecture 5.2]{DDR} and \cite[Conjecture 1.6]{DDS}, and which surprisingly states that {\bf the class $\mathcal{U}$ coincides with the class of all Abelian groups}. 

\medskip

We note that the result is preliminary announced in \cite{CD} as well as it is discussed and cited in \cite[Theorem 3.6]{S}. Specifically, we are prepared to prove the following main achievement.

\begin{theorem}\label{fine} Each group lies in $\mathcal{U}$, that is, any uniformly fully inert subgroup $H$ of a group $G$ is commensurable with some fully invariant subgroup of $G$.
\end{theorem}

\Pf It follows directly from a combination of Propositions~\ref{square} and \ref{2}, as expected.
\fine

The following comments are, hopefully, worthwhile.

\medskip

\noindent{\bf Remark.} We are able to give a quick sketch of a second, less conceptual, proof of Theorem~\ref{fine} with no exploiting the validity of Proposition~\ref{square}, but utilizing similar arguments to those as presented above combined with Proposition~\ref{2}.

\medskip

In fact, suppose $H$ is a uniformly fully inert subgroup in $G$. Then, by Fact 1 quoted above (and same as \cite[Theorem 4.6]{CDG}), we derive that $H\oplus H$ is a uniformly fully inert subgroup in $G\oplus G$, whence $$H\oplus H\sim C$$ for some characteristic subgroup $C$ in $G\oplus G$ (see \cite[Corollary 1.9]{DDR}). But Lemma~\ref{1} applies to get that $C$ is fully invariant in $G\oplus G$ and that it is also of the form $C=F\oplus F$, where $F$ is fully invariant in $G$. Indeed, we can write $$C=(C\cap G)\oplus (C\cap G)$$ and set $F=C\cap G$. So, $$H\oplus H\sim F\oplus F$$ and, in view of Fact 2 alluded to above, we infer at once that $H\sim F$, as wanted. The proof is completed.

\medskip

Nevertheless, it is worth mentioning that some recent work on that conjecture related to $p$-groups was done in \cite{Keef} as well.

\medskip
\medskip

\noindent{\bf Acknowledgement.} The authors would like to thank Prof. Patrick W. Keef from Whitman College, Walla Walla, WA, Unites States, for the professional checking of the truthfulness of the main result. They are also very thankful to the anonymous specialist referee for the very careful reading of the manuscript and the expert suggestions made, which significantly improved the presentation of the article and its main result.

\medskip
\medskip
\medskip
\medskip

\noindent{\bf Funding:} The scientific work of Andrey R. Chekhlov was supported by the Ministry of Science and Higher Education of Russia (agreement No. 075-02-2023-943). The scientific work of Peter V. Danchev was supported in part by the Bulgarian National Science Fund under Grant KP-06 No 32/1 of December 07, 2019, as well as by the Junta de Andaluc\'ia under Grant FQM 264, and by the BIDEB 2221 of T\"UB\'ITAK.

\end{document}